\documentclass{amsart}
\usepackage{amsfonts,amssymb,latexsym,amsthm,amsmath}

\usepackage{euscript}
\usepackage[arrow,matrix,curve]{xy}\SilentMatrices
\def\xyma{\xymatrix@M.7em}
\usepackage{float}

\usepackage{rotating}
\usepackage{verbatim}

\newcommand{\ilimit}{\,\varprojlim{}\!}
\newcommand{\dlimit}{\,\varinjlim{}\!}
\newcommand{\mono}{\rightarrowtail}

\usepackage{pb-diagram,pb-xy}

\newtheorem{theorem}{Theorem}

\begin{document}
\title{On transfinite nilpotence of the Vogel-Levine localization}
\author{Roman Mikhailov}
\address{Chebyshev Laboratory, St. Petersburg State University, 14th Line, 29b,
Saint Petersburg, 199178 Russia and St. Petersburg Department of
Steklov Mathematical Institute} \email{rmikhailov@mail.ru}
\thanks{This research is supported by the Chebyshev
Laboratory  (Department of Mathematics and Mechanics, St.
Petersburg State University)  under RF Government grant
11.G34.31.0026, and by JSC "Gazprom Neft", as well as by the RF
Presidential grant MD-381.2014.1.}
\begin{abstract}
We construct a finitely-presented group such that its Vogel-Levine
localization is not transfinitely nilpotent. This answers a
problem of J. P. Levine.
\end{abstract}
\maketitle

\section{Introduction}
In the series of papers \cite{Levine89},
\cite{Levine891},\cite{Levine1}, \cite{Levine2}, J. P. Levine
developed the theory of algebraic closure of groups and described
possible applications in geometric topology, as well as formulated
natural problems related to localizations and completions of
groups.

A group homomorphism $f: G\to H$ is called {\it 2-connected} if it
induces isomorphism on $H_1(-,\mathbb Z)$ and a surjection on
$H_2(-,\mathbb Z)$.  Denote by $\Omega$ the collection of all
2-connected homomorphisms $f: G\to H$ such that $G$ and $H$ are
finitely presented groups. This concept plays a fundamental role
in geometric topology, since a homology equivalence of connected
spaces induces a 2-connected homomorphism of their fundamental
groups (see \cite{CO1}, \cite{CO2} for geometric applications of
the theory of 2-connected homomorphisms).

 A group $\Gamma$ is {local} if given any diagram of
homomorphisms as follows, with $G \to H$ from $\Omega$,
\[
\begin{diagram}
\node{G} \arrow{e}\arrow{se} \node{H}\arrow{s,..}\\
\node[2]{\Gamma}
\end{diagram}
\]
there is a unique homomorphism $H \to \Gamma$ making the diagram
commute. The {\it Vogel-Levine localization} (also called an {\it
algebraic closure}) of a group $G$ is a group $L(G)$ endowed with
a homomorphism $l: G\to L(G)$, such that $L(G)$ is local and for
any local group $\Gamma$ and a homomorphism $f: G\to \Gamma$,
there is a unique homomorphism $p: L(G)\to \Gamma$ such that
$p\circ l=f$. The Vogel-Levine localization is an algebraic analog
of the localization of CW-complexes considered by J.-Y. Le Dimet
\cite{Di}.

The Vogel-Levine localization is a functor from all groups to the
local groups. The existence, uniqueness and different properties
of this functor are given in \cite{Levine1}, \cite{Levine89},
\cite{Levine891}. Recall some of the properties.

(i) Any homomorphism $f: G\to H$ from $\Omega$ induces an
isomorphism of localizations $L(G)\simeq L(H)$.

(ii) For any $G$, the localization $l: G\to L(G)$ is 2-connected.

For a finitely generated group $G$, the functor $L(G)$ lives in
the corner of the following square which plays a fundamental role
in the theory of localizations and completions of groups
\begin{equation}\label{square} \xyma{L(G) \ar@{->}[r] \ar@{->>}[d] &
E(G)\ar@{->>}[d]\\ \overline G \ar@{>->}[r] & \widehat G}
\end{equation}
Here $E(G)$ is the functor of HZ-localization defined by A. K.
Bousfield (see \cite{Bou}), $\widehat G:=\ilimit_i\ G/\gamma_i(G)$
is the functor of pro-nilpotent completion, $\overline G$ is the
$\omega$-closure (or {\it residually nilpotent algebraic closure},
see \cite{Levine89}, \cite{Levine1}, \cite{Levine2}).

Given a finitely presented group $G$, it is a difficult problem
how to describe the algebraic closure $L(G)$ and its
group-theoretical properties. One tool to recognize $L(G)$ is
given in \cite{CO2}. Suppose one can construct a sequence of
2-connected homomorphisms
$$
G\to K_1\to K_2\to \dots
$$
such that $K_i,\ i=1,2,\dots$ are finitely presented and $\dlimit\
K_i$ is a local group. Then $L(G)=\dlimit\ K_i$. This follows
immediately from the definition and uniqueness of the algebraic
closure. The seed of the idea to describe the group $L(G)$ as
injective limit of 2-connected maps is given in \cite{Di}. In
\cite{Di}, Le Dimet shows that the Vogel localization of a finite
CW-complex is a colimit of a countable sequence of finite
CW-complexes. The above tool to recognize $L(G)$ is an algebraic
analog of the construction from \cite{Di}.

For a group $G$, the lower central series are defined inductively
as follows: $$\gamma_1(G)=G,\
\gamma_{\alpha+1}(G)=[\gamma_\alpha(G),G]$$ and
$\gamma_\tau(G)=\cap_{\alpha<\tau}\gamma_\alpha(G)$ for a limit
ordinal $\tau$. A group $G$ is called {\it transfinitely
nilpotent} if $\gamma_\alpha(G)=1$ for some ordinal $\alpha$. The
vertical arrows in (\ref{square}) are quotients of $L(G)$ and
$E(G)$ by the intersections of (finite) lower central series
$\gamma_\omega:=\cap_i\gamma_i$.

For any group $G$, its HZ-localization $E(G)$ is transfinitely
nilpotent \cite{Bou}. In order to compare group-theoretical
properties of algebraic closures and HZ-localizations, J. P.
Levine asked the following (\cite{Levine1} (Problem 6 (b)): {\it
If $G$ is finitely-generated, is $L(G)$ transfinitely nilpotent?}
The following example answers this problem\footnote{In
\cite{Levine1}, J.P. Levine considered localization with respect
to 2-connected maps which are normally surjective. Observe that
for the group $H$ in theorem 1, $L(G)$ is the localization in that
sense as well, since all maps considered in the construction are
normally surjective}.

\begin{theorem}\label{th2}
Let $H=\langle a,b\ |\ a^{b^2}=aa^{3b}, [a,a^b]=1\rangle$. Then
$L(H)$ is not transfinitely nilpotent.
\end{theorem}

As usual, if $x,y,a_1,\dots,a_{k+1}$ are elements of a group $G$
we set $[x,y] =x^{-1}y^{-1}xy$, $x^y=y^{-1}xy$ and define
$$[a_1,a_2,\dots, a_{k+1}]=[[a_1,\dots,a_{k}], a_{k+1}]\ (\, k>1).$$

\section{Proof of theorem \ref{th2}}
The construction is based on the 1-relator group from \cite{Mikh}
$$G= \langle a,b\ |\ a^{b^2}=aa^{3b}\rangle
$$
with long lower central series. Our group $H$ is the quotient
$G/\gamma_\omega(G)$. The lower central series length of $G$ is
$\omega^2$. The $2\omega$-lower central quotient
$$G/\gamma_{2\omega}(G)=\langle a,b\ |\ a^{b^2}=aa^{3b},\
[a,a^b,a]=[a,a^b,a^b]=1\rangle$$ lives in the short exact sequence
$$
1\to \langle [a,a^b]\rangle \to G/\gamma_{2\omega}(G)\to H\to 1
$$
where $\langle [a,a^b]\rangle$ is the infinite cyclic group, $a\in
H$ acting on $[a,a^b]$ trivially and $b$ by inverting. The lower
central quotients $\gamma_{\omega+k}(G)/\gamma_{\omega+k+1}(G)$
are cyclic groups of order 2 for all $k=0,1,\dots$ with generators
$[a,a^b]^{2^k}.\gamma_{\omega+k+1}(G)$.

The proof consists of the following three steps:
\begin{enumerate}
\item Description of the $\omega$-closure $\overline H$ and the
proof that $H_2(\overline H)=0$. Since $\overline
H=L(H)/\gamma_\omega(L(H))$, the part of the 5-term sequence
$$
H_2(L(H))\to H_2(\overline H)\to
\gamma_{\omega}(L(H))/\gamma_{\omega+1}(L(H))\to 1
$$
implies that $\gamma_\omega(L(H))=\gamma_{\omega+1}(L(H))$. This
implies that $L(H)$ is residually nilpotent if and only if
$L(H)=\overline H$.\\
\item A construction of the sequence of finitely presented groups
$\Gamma_k',\ k=0,1,2,\dots$ and 2-connected homomorphisms
$$
H=\Gamma_0'\to \Gamma_1'\to \dots \to \Gamma_k'\to \Gamma_{k+1}'
\to \dots
$$
Denote the 2-connected maps from this sequence $h_k: H\to \Gamma_k',\ k=1,2,\dots$.\\
\item Proof that the limit $\Gamma:=\dlimit \Gamma_k'$ is a local
group. This implies that $L(H)=\Gamma$ and that the algebraic
closure $H\to L(H)$ is the limit map $$\dlimit\ h_k:\ H\to
\dlimit\ \Gamma_k'.$$ It will be shown that there is a natural
exact sequence $$ 1\to C_{2^{\infty}}\to L(H)\to \overline H\to 1
$$
That is, $\gamma_\omega(L(H))$ is non-trivial and isomorphic to
the 2-quasi-cyclic group $C_{2^\infty}:=\dlimit\{\mathbb
Z/2\mono\mathbb Z/4\mono \mathbb Z/8\mono\dots\}$.
\end{enumerate}

\vspace{.5cm}\noindent{\bf Step 1.} Let $N=\mathbb Z\oplus \mathbb
Z$ be the $\mathbb Z[\langle b\rangle]$-module generated by $a$
and $a^b$.  The generator $b$ of the cyclic quotient of $H$ acts
on $N$ as the matrix
\begin{equation}
U:=\begin{pmatrix}
0 & 1 \\
1 & 3
\end{pmatrix}
\end{equation}
Recall that (see prop. 3.2. from \cite{Levine2}, and {\it
telescope theorem} from \cite{BMO1}, \cite{BMO2}) the
$\omega$-closure $\overline H$ of $H$ has the following natural
description
$$
\overline H=N_S\rtimes \mathbb \langle b\rangle,
$$
where $N_S$ is the $S$-localization of $N$ with $S:=1+\Delta,\
\Delta=\ker\{\mathbb Z[\langle b\rangle]\to \mathbb Z\}.$ By
construction, $N_S$ is the direct limit
$$
N_S=\dlimit\{N\buildrel{s_1}\over\longrightarrow
N\buildrel{s_2}\over\longrightarrow \dots \}
$$
where $\{s_1,s_2,\dots\}$ covers all elements from $S$.

For an element $s\in S$, the $s$-map $N\buildrel{s}\over\to N$ is
injective. This follows from the residual nilpotence of $H$ (see
\cite{Mikh} for the proof that $H$ is residually nilpotent).
Indeed, if $s(n)=0,$ for some $n\in N, n\geq 0,$ then $n\in
\gamma_\omega(H)=1$.

The multiplication by $s\in S$, induces a map of exterior squares
$$
(\mathbb Z\simeq \langle a\wedge a^b\rangle \simeq)\Lambda^2(N)\to
\Lambda^2(N).
$$
Any homomorphism $\mathbb Z\to \mathbb Z$ is a muliplication with
some number. In our case denote this number by $|s|$. Let
$|s|=2^{p(s)}v$ for an odd $v$. Observe that $|s|\neq 0$, since,
as observed above, any $s$-map $N\buildrel{s}\over\to N$ is
injective.

Since $H_2(H)=\mathbb Z/2$, the map $s$ induces a zero map
$H_2(H)\to H_2(H)$ if and only if $p(s)\neq 0.$ The simplest
example of an element from $S$, with zero induced map
$H_2(H)\buildrel{s_*}\over\rightarrow H_2(H)$ is
$$
s=1-b+b^2.
$$
Consider an element from $S$ of the form
$$
s=n_0+n_1b+\dots+n_lb^l, \ n_0+\dots+n_l=1.
$$
For $a\in N$, we will use the natural notation
$$
a^s:=a^{n_0}a^{n_1b}\dots a^{n_lb^l}.
$$
A technical exercise is to show that
$$
|s|\equiv 1+\sum_{j>i,\ 3\nmid (j-i)} n_in_j \mod 2
$$
Since there are infinitely many indecomposed elements in $S$,
which induce the zero map on $H_2(H)$, we conclude that
$$
H_2(\overline H)=\dlimit_{s_1,s_2,\dots} H_2(H)=0.
$$
Using the general relation between Vogel-Levine localization and
$\omega$-closure, $\overline H=L(H)/\gamma_{\omega}(L(H)),$ the
5-term sequence implies that
$$
\gamma_\omega(L(H))=\gamma_{\omega+1}(L(H)).
$$

\vspace{.5cm}\noindent{\bf Step 2.} Denote the following groups
$$
\Gamma_k:=\langle a,b\ |\ a^{b^2}=aa^{3b},
[a,a^b,a]=[a,a^b,a^b]=1, [a,a^b,\underbrace {b,\dots,
b}_k]=1\rangle,\ k\geq 1
$$

Consider an element $s\in S$. For a given $k\geq 0$, we will
construct a homomorphism $\phi_s: \Gamma_k\to \Gamma_{k+p(s)}$
such that there is a commutative diagram
$$
\xyma{\gamma_\omega(\Gamma_k)(\simeq \mathbb Z/2^k)\ar@{>->}[r]
\ar@{>->}[d] &
\gamma_\omega(\Gamma_{k+p(s)})(\simeq \mathbb Z/2^{k+p(s)})\ar@{>->}[d] \\ \Gamma_k \ar@{->}[r]^{\phi_s}\ar@{->>}[d] & \Gamma_{k+p(s)}\ar@{->>}[d]\\
H \ar@{->}[r]^{s_*} & H}
$$
and $\phi_s$ induces an isomorphism $H_2(\Gamma_k)\to
H_2(\Gamma_{k+p(s)}).$

Since the map $N\buildrel{s}\over\to N$ induces a well-defined
homomorphism $(N\rtimes\langle b\rangle=)H\to H$, there exists
$l\geq 0$, such that
\begin{equation}\label{rl}
(a^s)^{b^2}=a^s(a^{3s})^b[a,a^b]^l \end{equation} in $\Gamma_k$.
We define the homomorphism $\phi_s: \Gamma_k\to \Gamma_{k+p(s)}$
as
\begin{align*}
& a\mapsto a^s[a,a^b]^r\\
& b\mapsto b
\end{align*}
with $r$ such that $3r\equiv l\mod 2^k$. Lets check that $\phi_s$
is well-defined. Indeed, the relation (\ref{rl}) implies that, in
$\Gamma_{k+p(s)}$,
$$
(a^s[a,a^b]^r)^{b^2}=(a^s[a,a^b]^r)(a^s[a,a^b]^r)^{3b}
$$
The group $\Gamma_k$ has another relation which we have to check.
We have
$$
[[a,a^b]^s,\underbrace {b,\dots,
b}_k]=([a,a^b]^s)^{2^k}=[a,a^b]^{|s|2^k}=1
$$
in $\Gamma_{k+p(s)}$. The relations $[a,a^b,a]=[a,a^b,a^b]=1$ are
preserved by the considered map. Thus, the homomorphism $\phi_s$
is well-defined. Observe that the homomorphism $\phi_s$ is
normally surjective, i.e. the normal closure of the image of
$\phi_s$ equals to $\Gamma_{k+p(s)}$.

 The homology group $H_2(\Gamma_k)$ is isomorphic to
$\mathbb Z/2$. Looking at presentation of the second homology
$H_2(\Gamma_k)$ via the Hopf formula $H_2(\Gamma_k)=\frac{R\cap
[F,F]}{[R,F]},$ with $F=F(a,b)$, we describe the generator of
$H_2(\Gamma_k)$ as the coset $$[a,a^b,\underbrace {b,\dots,
b}_k].[R,F]=[a,a^b]^{2^k}.[R,F].$$The image of the generator of
$H_2(\Gamma_k)$ under the map induced by $\phi_s$ is non-trivial
in $H_2(\Gamma_{k+p(s)})$, hence the induced map $$
H_2(\Gamma_k)(\simeq \mathbb Z/2)\to H_2(\Gamma_{k+p(s)})(\simeq
\mathbb Z/2)
$$
is an isomorphism.

Lets illustrate the above construction for the particular case
$s=1-b+b^2$. In this case, $|s|=12$ and the induced map $H_2(H)\to
H_2(H)$ is zero. We claim that there is a commutative diagram
$$
\xyma{& \Gamma_2\ar@{->>}[d]\\ H\ar@{->}[r]^{s_*}
\ar@{->}[ur]^\phi & H}
$$
such that $\phi$ induces isomorphism $H_2(H)\to H_2(\Gamma_2)$ and
the vertical map is $\Gamma_2\to
\Gamma_2/\gamma_{\omega}(\Gamma_2)=H$. The map $\phi$ is defined
as
\begin{align*}
& a\mapsto aa^{-b}a^{b^2}[a,a^b]\\
& b\mapsto b
\end{align*}
One can easily see that, in $\Gamma_2$,
$$
(aa^{-b}a^{b^2}[a,a^b])^{b^2}=aa^{-b}a^{b^2}[a,a^b](aa^{-b}a^{b^2}[a,a^b])^{3b}.
$$
hence the map $\phi$ is well-defined.

To finalize the Step 2, we conclude that, for a sequence of
elements $\{s_1,s_2,\dots\}$, there exists an infinite tower
$$
\xyma{& \mathbb Z/2^{|p(s_1)|} \ar@{>->}[d] \ar@{>->}[r] & \mathbb
Z/2^{|p(s_1)|+|p(s_2)|} \ar@{>->}[d] \ar@{>->}[r] & \dots\\
\Gamma_0\ar@{->}[r]^{\phi_{s_1}}\ar@{->}[d]^{=} & \Gamma_{|p(s_1)|} \ar@{->}[r]^{\phi_{s_2}}\ar@{->>}[d] & \Gamma_{|p(s_1)|+|p(s_2)|}\ar@{->}[d] \ar@{->}[r] & \dots\\
H \ar@{->}[r]^{s_1} & H \ar@{->}[r]^{s_2} & H \ar@{->}[r] & \dots
}
$$
All the homomorphisms $\phi_{s_i}$ are 2-connected. The group
$\Gamma:=\dlimit_{\{s_1,s_2,\dots\}}\Gamma_{s_i}$, lies in the
short exact sequence
$$ 1\to C_{2^{\infty}}\to \Gamma\to \overline H\to 1
$$
and the action of $\overline H=N_S\rtimes \langle b\rangle$ on the
quasi-cyclic group $C_{2^{\infty}}$ is given as follows:
$$
y\circ b=-y,\ y\circ n=y,\ n\in N_S.
$$

\vspace{.5cm}\noindent{\bf Step 3.} In order to show that $\Gamma$
is local, recall the definition of local Cohn modules. For a group
$Q$, let $M$ be a $\mathbb Z[Q]$-module. We call $M$ a {\it local
Cohn module} if, for every map $t: F_1\to F_2$ of finitely
generated free $\mathbb Z[Q]$-modules of the same rank, such that
the induced map $1\otimes_{\mathbb Z[Q]}\mathbb Z:
F_1\otimes_{\mathbb Z[Q]}\mathbb Z\to F_2\otimes_{\mathbb
Z[Q]}\mathbb Z$ is an isomorphism, and a morphism of $\mathbb
Z[Q]$-modules $\alpha: F_1\to M,$ there is a unique morphism
$\beta: F_2\to M$ such that $\beta\circ t=\alpha.$

Recall the following result (see \cite{CO1}), \cite{CO2}). Let $Q$
be a local group and $M$ a Cohn local $\mathbb Z[Q]$-module. For
any extension
$$
0\to M\to \widetilde Q\to Q\to 1,
$$
the group $\widetilde Q$ is local.

Observe that $C_{2^\infty}$ is the direct limit of $\mathbb
Z[\overline H]$-modules $$\dlimit\{\mathbb Z/2\mono\mathbb
Z/4\mono \mathbb Z/8\mono\dots\}$$ Every submodule $\mathbb
Z/2^k,\ k\geq 1$ is nilpotent $\mathbb Z[\overline H]$-module, $
(\mathbb Z/2^k)\Delta^k(\overline H)=0. $ Every nilpotent module
is Cohn local and since $C_{2^\infty}$ is a direct limit of
nilpotent modules, we conclude that $C_{2^\infty}$ is Cohn local
module. Hence $\Gamma$ is a local group and the Step 3 is
complete. This completes the proof of theorem.

\vspace{.5cm}\noindent {\it Acknowledgements.} The author thanks
S.O. Ivanov and K. Orr for discussions related to the subject of
the paper.

\end{document}